\documentclass[11pt]{article}
\usepackage{graphicx}
\usepackage{amsmath,amsfonts,amssymb,bm}

\newcommand{\CA}{\mathcal C}
\newcommand{\RCA}{\mathcal R}

\newcommand{\Eq}[1]{Eq.~(\ref{#1})}
\newcommand{\Fig}[1]{Fig.~\ref{fig:#1}}

\newcommand{\swap}{\mathsf X}

\title{On number of nonzero cells in some
2D reversible second-order cellular automata}
\date{\today}
\author{Alexander Yu.\ Vlasov}
\begin{document}
\sloppy
\maketitle
\begin{abstract}
Recursive equations for the number of cells with nonzero values at
$n$-th step for some two-dimensional reversible second-order cellular 
automata are proved in this work. Initial configuration is a single 
cell with the value one and all others zero. 
\end{abstract}

\section{Introduction}

Any cellular automaton (CA) with two states and 
local transition rule \mbox{$c \mapsto f[c]$} can be used
for definition of a reversible second-order CA
with new rule $F$ on a pair \cite{TM90,WCA}
\begin{equation}
 F : (c,c') \mapsto (f[c] + c' \mod 2\, , c).
 \label{2nd}
\end{equation}
An inverse rule is 
\begin{equation}
 F^{-1} : (c',c) \mapsto (c,f[c] + c' \mod 2).
 \label{inv2nd}
\end{equation}
and also may be rewritten
\begin{equation}
 F^{-1}  = \swap F \swap ,
 \label{invswap}
\end{equation}
where $\swap$ is exchange of states
\begin{equation}
 \swap : (c,c') \mapsto (c', c).
 \label{swap}
\end{equation}
The CA acting on pairs of binary states $(c,c')$ can be considered as
four-state CA due to simple correspondence $(c,c') \mapsto c+2c'$.

Let us denote $c_{i,j}$ state of a cell with notations 
\begin{eqnarray*}
\Sigma^\times c_{i,j} &=& c_{i-1,j-1} + c_{i-1,j+1} +  c_{i+1,j-1} + c_{i+1,j+1}, \\
\Sigma^+ c_{i,j} &=& c_{i,j-1} + c_{i,j+1} +  c_{i+1,j} + c_{i-1,j}.
\end{eqnarray*}
Let us consider a few different two-dimensional CA
\begin{enumerate}
\item $\CA_1$ with local rule:
$c_{i,j} \mapsto \Sigma^\times c_{i,j} \mod 2$
\item $\CA_2$ with local rule:
$c_{i,j} \mapsto \Sigma^+ c_{i,j} \mod 2$:
\item $\CA_3$ with local rule:
$
c_{i,j} \mapsto \begin{cases}
1, & \Sigma^+ c_{i,j} = 1\\
0, & \text{otherwise}
\end{cases}
$

$\CA_3'$ with local rule:
$
c_{i,j} \mapsto \begin{cases}
1, & \Sigma^+ c_{i,j} = 1 \text{ and } \Sigma^\times c_{i,j} = 0 \\
0, & \text{otherwise}
\end{cases}
$
\end{enumerate}
and second-order reversible CA
$\RCA_1, \RCA_2, \RCA_3, \RCA_3'$ 
derived from them using \Eq{2nd}. 

If to start with a single cell with value one and all others zero,
then total number of cells with nonzero values at $n$-th stage
is some sequence $R(n)$. It is also possible to consider sequences 
$R_k(n)$, $k=1,2,3$ for number of cells with value $k$. 

The sequence was initially introduced due to consideration
of ``noise'' in computationally universal CA $\RCA_3'$, but it is 
shown below, that for other three CA the sequences are the same 
and $R_3(n) = 0$. Due to definition of second-order CA \Eq{2nd} 
a simple property is true
\begin{equation}
R_2(n+1)=R_1(n)
\label{R1R2}
\end{equation} and so 
\begin{equation}
R(n)=R_2(n)+R_2(n+1).
\label{Rsum}
\end{equation}
Initial terms of the sequences are represented in the table below:
\begin{equation}
\setlength{\arraycolsep}{0.3em}
\begin{array}{|l|r|r|r|r|r|r|r|r|r|r|r|r|r|r|r|r|}\hline
  n & 0 & 1 & 2 & 3  & 4  & 5  & 6  & 7  & 8  & 9  & 10 & 11 & 12 & 13 & 14 & 15\\ \hline
  R & 1 & 5 & 9 & 21 & 25 & 29 & 41 & 85 & 89 & 61 & 65 & 109& 121& 125& 169& 341\\ 
R_1 & 1 & 4 & 5 & 16 & 9  & 20 & 21 & 64 & 25 & 36 & 29 & 80 & 41 & 84 & 85 & 256\\ 
R_2 & 0 & 1 & 4 & 5  & 16 & 9  & 20 & 21 & 64 & 25 & 36 & 29 & 80 & 41 & 84 & 85\\ \hline
\end{array}
\label{tabR}
\end{equation}

\section{Recursive equations for numbers of cells}

Few recursive equations are proved in this paper:
\begin{equation}
 R(0)=1, \quad R(2^k+j) = 4 R(j) + R(2^k-j-1),\quad 0 \le j < 2^k.
\label{recR} 
\end{equation}
\begin{equation}
 R_1(-1)=0, \quad R_1(0)=1, 
 \quad R_1(2^k+j) = 4 R_1(j) + R_1(2^k-j-2).
\label{recR1} 
\end{equation}
The negative value of $n$ can be used because CA are reversible.
Due to \Eq{R1R2} last formula is equivalent with
\begin{equation}
 R_2(0)=0, \quad R_2(1)=1, 
 \quad R_2(2^k+j) = 4 R_2(j) + R_2(2^k-j).
\label{recR2} 
\end{equation}
Both \Eq{recR1} and \Eq{recR} are simply derived from the equation \Eq{recR2}:
\begin{eqnarray*}
 R_1(2^k+j) &=& R_2(2^k+j+1) = 4 R_2(j+1) + R_2(2^k-j-1)\\ 
            &=& 4 R_1(j) + R_1(2^k-j-2), 
\end{eqnarray*}
\begin{eqnarray*}
 R(2^k+j) &=& R_2(2^k+j)+R_2(2^k+j+1)\\ 
          &=& 4 \bigl(R_2(j) + R_2(j+1)\bigr) + R_2(2^k-j)+R_2(2^k-j-1)\\ 
          &=& 4 R(j) + R(2^k-j-1). 
\end{eqnarray*}

An {\em alternative form of recursive equations} is also valid for $R_1$ and $R_2$:
\begin{equation}
 R_1(2 n + 1) = 4 R_1(n), \quad R_1(2 n + 2) =  R_1(n) + R_1(n+1), 
\label{rec2R1} 
\end{equation}
\begin{equation}
 R_2(2 n) = 4 R_2(n), \quad R_2(2 n + 1) =  R_2(n) + R_2(n+1). 
\label{rec2R2} 
\end{equation}
These equations are equivalent due to \Eq{R1R2} and together with \Eq{Rsum} 
imply a simple relation between the sequences
\begin{equation}
 R(n) = R_1(2n) = R_2(2n+1).
 \label{Rdub}
\end{equation}

Equations \Eq{rec2R1} and \Eq{rec2R2} can be proved {\em by induction}
using \Eq{recR1} and \Eq{recR2} respectively. Due to \Eq{R1R2}, 
it is enough to consider only one of them.

The \Eq{recR2} holds for $0\le k<4$. Assume \Eq{recR2} holds
for any $n < k$, $k = 2^m+j$ with $m >0 $, $0 < j \le 2^m$.
\Eq{recR2} allows us to express $R_2$ as a linear combination
with terms smaller than $k$ and to show that the equation holds also for 
$n = k$:
\begin{eqnarray*}
R_2(2n) &=& R_2(2^{m+1}+2j) = 4R_2(2j)+ R_2(2^{m+1}-2j)\\
         &=& 4(4R_2(j)+ R_2(2^m-j)) = 4 R_2(2^m+j) = 4 R_2(n),
\end{eqnarray*}
\begin{eqnarray*}         
R_2(2n+1) &=& R_2(2^{m+1}+2j+1) = 4R_2(2j+1)+ R_2(2^{m+1}-2j-1)\\ 
          &=& 4 R_2(j)+ 4 R_2(j+1)+ R_2(2^m-j)+  R_2(2^m-j-1)\\     
          &=& R_2(2^m+j)+R_2(2^m+j+1) = R_2(n)+R_2(n+1),
\end{eqnarray*}
where $j < j+1 = k-2^m+1 < k$, $2^m-j-1 < 2^m-j = k-2j < k$. $\square$ 

\smallskip

It remains to prove \Eq{recR2}. 
The recursion is proved below for simpler case with CA 
$\RCA_1$ and $\RCA_2$  
with straightforward demonstration of equivalence for CA
$\RCA_3$ and $\RCA_3'$.

\section{Properties of initial two-state CA} 

Let us start with consideration of $\CA_1$ and
$\CA_2$. These CA are linear (additive) \cite{alinCA,statCA}, {\em i.e.}
for any two configurations $a$ and $b$ local rule defines 
global map $f$ with property
\begin{equation}
 f(a \oplus b) = f(a) \oplus f(b),
 \label{lin2}
\end{equation}
where $a \oplus b = a \bigtriangleup b =
(a \cup b)\smallsetminus(a \cap b) =
(a \smallsetminus b)\cup (a \smallsetminus b)$ is symmetric difference
configurations $a$ and $b$ considered as sets (regions) of cells with
unit values.

A configuration of 2D CA can be described with (characteristic) polynomial
\begin{equation}
 p[c] \equiv p_{x,y}[c] = \sum_{i,j=-\infty}^\infty c_{i,j} x^i y^j
 \label{poly}
\end{equation}
and \Eq{lin2} corresponds to
\begin{equation}
 p[f(a \oplus b)] = p[f(a)] \oplus p[f(b)] 
 \equiv p[f(a)] + p[f(b)] \mod 2.
 \label{polin2}
\end{equation}
It is convenient further for CA with two states
to treat \Eq{poly} as a polynomial over $GF(2)$.

Let us consider evolution of pattern $\Delta_{0,0}$ with 
single nonzero cell $c_{0,0}=1$ for CA $\CA_1$. 
It can be described using equation for global transition rule
\begin{eqnarray}
 \CA_1: p_{x,y}[c] &\mapsto &
  (x^{-1}y^{-1} + x y^{-1} + x^{-1}y + xy)\, p_{x,y}[c] \notag\\
 & = & (x^{-1}+x )(y^{-1} + y)\,P_{x,y}[c].
\label{C1poly} 
\end{eqnarray}
Here treatment of $p_{x,y}[c]$ as a polynomial over $GF(2)$
is especially useful and after $n$ 
steps due to \Eq{C1poly}
\begin{equation}
  p_{x,y}[\CA_1^n(c)] = (x^{-1}+x)^n (y^{-1} + y)^n\,p_{x,y}[c].
 \label{C1poln}  
\end{equation}
The polynomial of pattern $\Delta_{0,0}$ is $p_{x,y}[c]=p[\Delta_{0,0}]=1$
and the \Eq{C1poln} corresponds to decomposition $p_{x,y}=p_x p_y$ on two
characteristic polynomials $p_x=x^{-1}+x$ of 1D 
cellular automata with local rule \cite{statCA,algprop} 
\begin{eqnarray}
a_i \mapsto a_{i-1}+a_{i+1} \mod 2
\label{rule90}
\end{eqnarray}
also known as ``rule 90'' \cite{statCA}
and initial pattern $\Delta_0$ with single nonzero cell $a_0=1$.
The number of cells on $k$-th step may be described by equation
\begin{equation}
 \mathcal N_k = 2^{\ell(k)},
\label{N90} 
\end{equation}
there $\ell(k)$ is number of units in binary decomposition of $k$ \cite{statCA}.

The polynomial is over $GF(2)$ and a property used further  
\begin{equation}
(x^{-1}+x)^{2^k} = (x^{-2^k}+x^{2^k}) \mod 2
\label{binom2k} 
\end{equation}
is simply derived using recursion on $k$: 
$$(x^{-1}+x)^{2^{k+1}} = (x^{-2^k}+x^{2^k})^2.$$

\Eq{binom2k} can be used for inductive proof of \Eq{N90}. 
For $k=0$ \Eq{N90} holds: $\mathcal N_k=1$. 
Assume $N_k = 2^{\ell(k)}$ for $k=0,\ldots,2^n$. 
For $k' = k+2^n$ characteristic polynomial is 
$p_{k'}(x) = (x^{-2^n}+x^{2^n})p_k(x)$ and because
$x^{-2^n}p_k(x)$ and $x^{2^n}p_k(x)$ are not ``overlapped,''
$N_{k'} = 2 N_k$. Due to $\ell(k') = \ell(k)+1$
for $k'+2^n$: $N_{k'} = 2 N_k = 2^{\ell(k)+1}$. 
So, \Eq{N90} holds for $k=0,\ldots,2^{n+1}$. $\square$

The decomposition \Eq{C1poly} produces some simplification with 
comparison to $\CA_2$
\begin{equation}
 \CA_2: p_{x,y}[c] \mapsto 
  (x^{-1} + x + y^{-1} + y)\, p_{x,y}[c] .
\label{C2poly} 
\end{equation}

On the other hand, $\CA_1$ ($\RCA_1$) may be considered as two independent copies
of $\CA_2$ ($\RCA_2$) on two {\em ``diagonal'' sublattices} corresponding $c_{i,j}$ with even
and odd $i+j$ respectively: 
\begin{equation}
c'_{i,j} = c_{i+j,i-j},\quad c''_{i,j} = c_{i+j+1,i-j}.
\label{C2C1}
\end{equation}
Visually, they correspond to cells with black and white colors on checkerboard pattern 
after $\pi/4$ rotation of the board.

Because $c_{0,0}$ belongs to even sublattice $c'$, configuration of $\CA_1$ after any
$n$ steps always belongs to $c'$ and it is equivalent with $\CA_2$ acting
on the diagonal sublattice.

\smallskip

Due to \Eq{binom2k} and \Eq{C1poln} application of $2^k$ steps of $\CA_1$
to arbitrary configuration $c$ may be expressed as
\begin{equation}
\CA_1^{2^k}: p_{x,y}[c] \mapsto 
 (x^{-2^k}y^{-2^k}+x^{2^k}y^{-2^k} + x^{-2^k}y^{2^k}+x^{2^k}y^{2^k})\,p_{x,y}[c]
\label{C1expmov}
\end{equation}
and analogue property can be proved for $\CA_2$
\begin{equation}
\CA_2^{2^k}: p_{x,y}[c] \mapsto (x^{-2^k}+x^{2^k}+y^{-2^k} + y^{2^k})\,p_{x,y}[c].
\label{C2expmov}
\end{equation}
So patterns bounded by $2^k \times 2^k$ are replicated into four copies after $2^k$
steps both for $\CA_1$ and $\CA_2$. For $\CA_1$ coordinates of four copies are shifted
due to \Eq{C1expmov} as $(-2^k,-2^k)$, $(-2^k,+2^k)$, $(+2^k,-2^k)$, $(+2^k,+2^k)$ and
for $\CA_2$ due to \Eq{C2expmov} the shifts are 
$(-2^k,0)$, $(0,-2^k)$, $(+2^k,0)$, $(0,+2^k)$.
Such CA with replicating property was initially considered by E. Fredkin 
in 1970s \cite{CAPhy}.

For $\CA_1$ and $\CA_2$ an analogue of \Eq{N90} is true 
\begin{equation}
 \mathcal N_k = 4^{\ell(k)}.
\label{N90sqr} 
\end{equation}

Configuration of $\CA_1$ is represented as product 
\Eq{C1poln} $p(x,y)=p(x)\,p(y)$ of two
``rule 90'' CA and \Eq{N90sqr} can be derived directly from \Eq{N90}
\mbox{$\mathcal N_k = (2^{\ell(k)})^2$}. 
In more general case for such products of two 1D 
configurations $a$ and $b$ an equation 
$\mathcal N(a \cdot b) = \mathcal N(a) \mathcal N(b)$ can be used,
where $c = a_x \cdot b_y$ is 2D configuration with values
of cells $c_{i,j} = a_i b_j$.

\smallskip

A direct proof by induction for $\CA_2$ or $\CA_1$ is also
useful due to similarity with further approach to 
second-order CA. 
For $k=0$ \Eq{N90sqr} holds: \mbox{$\mathcal N_k=1$}. 
Assume $\mathcal N_k = 4^{\ell(k)}$ for $k=0,\ldots,2^n$. 
For $k' = k+2^n$ characteristic polynomial for $\CA_2$ 
satisfies \Eq{C2expmov} 
$$p_{k'}(x,y) = (x^{-2^n}+x^{2^n}+y^{-2^n}+y^{2^n})\,p_k(x,y),$$
and describes four shifted nonoverlapping copies of region $P_k(x,y)$.
So, $\mathcal N_{k'} = 4 \mathcal N_k = 4^{\ell(k)+1} = 4^{\ell(k')}$
and \Eq{N90sqr} holds for $k=0,\ldots,2^{n+1}$. $\square$

Similar proof by induction for $\CA_1$ uses \Eq{C1expmov}.

\section{Evolution of derived second-order CA}

A second-order CA corresponds to pair of polynomials
$\bigl(p_1(x,y),p_2(x,y)\bigr)$. For second-order CA 
derived from CA with two states described
by polynomials over GF(2) local rule \Eq{2nd}    
can be simply rewritten as a global one 
\begin{equation}
 \bigl(p_1(x,y)\, ,p_2(x,y)\bigr) \mapsto 
 \bigl(f[p_1(x,y)] + p_2(x,y)\, , p_2(x,y)\bigr).
 \label{plin2nd}
\end{equation}

For $\CA_1$, $\CA_2$ due to \Eq{C1poly} and \Eq{C2poly}
\begin{equation}
f[p(x,y)] = T(x,y) p(x,y) \mod 2
\label{Tp}
\end{equation} 
with 
\begin{equation}
T_{\CA_1}(x,y) = (x^{-1}+x )(y^{-1} + y),
\label{TC1}
\end{equation}
\begin{equation}
T_{\CA_2}(x,y) = (x^{-1}+x +y^{-1} + y).
\label{TC2}
\end{equation}

Let us prove that for  $\CA_1$, $\CA_2$ with
initial configuration $C_0 = \Delta_{0,0}$ 
with single nonempty cell $c_{0,0} = (1,0) \equiv \bf 1$  
after $k$ steps the configuration is described 
by polynomial
\begin{equation}
  P[C_k] = \bigl(f_{k+1}(T),f_k(T)\bigr),
  \label{PFT}
\end{equation}
where $f_k(t)$ are polynomials over GF(2) defined
using recursive equation
\begin{equation}
  f_{k+1}(t) = t f_k(t) + f_{k-1}(t),
  \quad  f_0 = 0, \quad f_1 = 1
\label{fibpol}  
\end{equation}
and $f_k(T)$ is application of the polynomial to
$T(x,y)$ \Eq{Tp} also considered over GF(2).
For $k \le 1$ \Eq{PFT} holds $P[\Delta_{0,0}] = (1,0)$
Assume \Eq{PFT} holds for $0,\ldots,k$, for $k+1$
due to \Eq{plin2nd}
$$
\bigl(f_k(T),f_{k - 1}(T)\bigr) \mapsto
\bigl(T\,f_k(T) + f_{k - 1}(T),f_k(T)\bigr) =
\bigl(f_{k+1}(T),f_k(T)\bigr). ~\square
$$

The \Eq{fibpol} defines {\em Fibonacci polynomials}. 
The {\em Lucas polynomials} (also used below) are
defined by the same recursive equation with other
initial conditions \cite{FibLuc}
\begin{equation}
  l_{k+1}(t) = t l_k(t) + l_{k-1}(t),
  \quad  l_0 = 2, \quad l_1 = t,
\label{lucpol}  
\end{equation}
\begin{equation}
l_k(t) = f_{k+1}(t) + f_{k-1}(t)= t f_k(t) + 2 f_{k-1}(t) 
\label{fibluc}
\end{equation}
with simpler correspondence over GF(2) 
\begin{equation}
t f_k(t) = l_k(t) \mod 2.
\label{fibluc2}
\end{equation}
Some relations with Lucas and Fibonacci polynomials \cite{FibLuc}
are useful further
\begin{equation}
 f_{m+n}(t) = f_m(t)l_n(t)+(-1)^{n+1} f_{m-n}(t),
 \label{fiblucrec}
\end{equation}
\begin{equation}
 f_{m+n+1}(t) = f_{m+1}(t)f_{n+1}(t)+ f_m(t)f_n(t).
 \label{fib2rec}
\end{equation}

For GF(2) multiplier $(-1)^{n+1}$ can be omitted and due
to relation \Eq{fibluc2} from \Eq{fiblucrec} for polynomials 
$f_k(T)$ over GF(2) follows 
\begin{equation}
f_{m+n}(T) = T f_m(T)f_n(T) + f_{m-n}(T).
\label{TFTrec}
\end{equation}
For $m=n$ \Eq{TFTrec} gives
\begin{equation}
f_{2n}(T) = T\,f_n^2(T)
\label{TFT2n}
\end{equation}
and \Eq{fib2rec} gives for $m=n$ 
\begin{equation}
 f_{2n+1}(t) = f_{n+1}^2(t) + f_n^2(t).
 \label{fib2n1}
\end{equation}
It again may be modified for polynomials over GF(2) 
\begin{equation}
 f_{2n+1}(T) = \bigl(f_{n+1}(T) + f_n(T)\bigr)^2.
 \label{FFT2n1}
\end{equation}

\smallskip

Let us show for polynomials over GF(2)
\begin{equation}
  f_{2^k}(T) = T^{2^k-1}.
  \label{Fpow2}
\end{equation}
It holds for $k=0$ and for $k+1$ due to \Eq{TFT2n}
$$
 f_{2^{k+1}}(T) = T\,f_{2^k}^2(T) = T^{2(2^k-1)+1} = T^{2^{k+1}-1}. ~\square
$$

Let us consider $f_{2^k+j}$ with $ j < 2^k$. 
Due to \Eq{TFTrec} and \Eq{Fpow2}
\begin{equation}
f_{2^k+j}(T) = T f_{2^k}(T) f_j(T) + f_{2^k-j}(T) = T^{2^k} f_j(T) + f_{2^k-j}(T).
\label{f2kj}
\end{equation}

\section{Proof of recursive equations}

A state of cell in the second-order CA for pair $(b_1,b_2)$ 
was encoded as $b_1 + 2 b_2$. The values one and two correspond to pairs
$(1,0)$ and $(0,1)$ respectively. 

Let us discuss distribution of cells with different values and  show, 
that $R_3(n)=0$, {\em i.e.} pair $(1,1)$ never appears for initial configuration with
single cell $(1,0)$.

The simpler way is to consider $\RCA_2$ ($\CA_2$) with checkerboard 
coloring already used earlier. Consider configuration with 
properties: \label{colprop}
\begin{enumerate}
  \item cells may not have state $(1,1)$
  \item all cells with the same state have the same color
\end{enumerate}
Show that these properties are valid after next step.
Let us denote $c_1$ and $c_2$ configurations
corresponding to set of cells with nonzero 
first and second elements of pair $(b_1,b_2)$ respectively
The properties above claim that configurations $c_1$ and
$c_2$ belong to {\em diagonal sublattices} with opposite colors.

The sublattices are represented by polynomials with odd and even 
degrees, so configurations with properties above correspond
to either (even,odd) or (odd,even) pairs of polynomials. 
The operator $T_{\CA_2}$ \Eq{TC2} changes degree of monomial
on unit and so \Eq{Tp} exchanges odd and even polynomials and
\Eq{plin2nd} maps configuration (odd,even) into (even,odd)
and vise versa. Initial configuration also has desired 
properties and so equation $R_3(k)=0$ is proved by induction. $\square$

It is more convenient sometimes to use $\CA_1$ instead of
$\CA_2$ and it is possible to introduce analogues of structures
discussed below. It was already mentioned that $\CA_2$ itself
corresponds to {\em diagonal sublattice} of $\CA_1$ and 
so notion of cells with the {\em ``same color''} needs for
some clarification.

\begin{figure}[htb]
\begin{center}
\parbox[t]{0.4\textwidth}{
\includegraphics[scale=0.5]{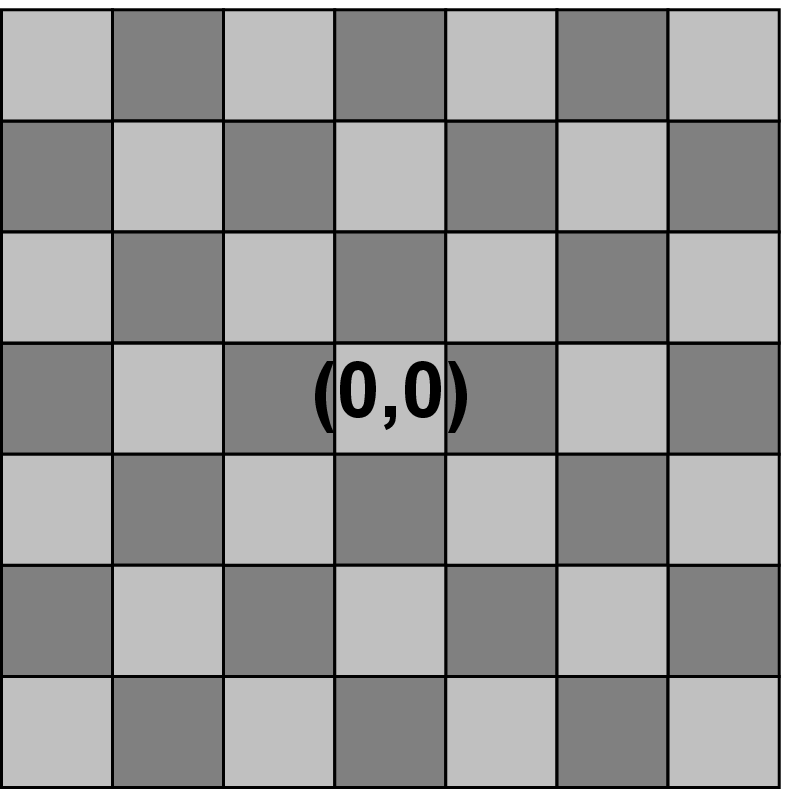} $\CA_2$}\hfil
\parbox[t]{0.4\textwidth}{
\includegraphics[scale=0.5]{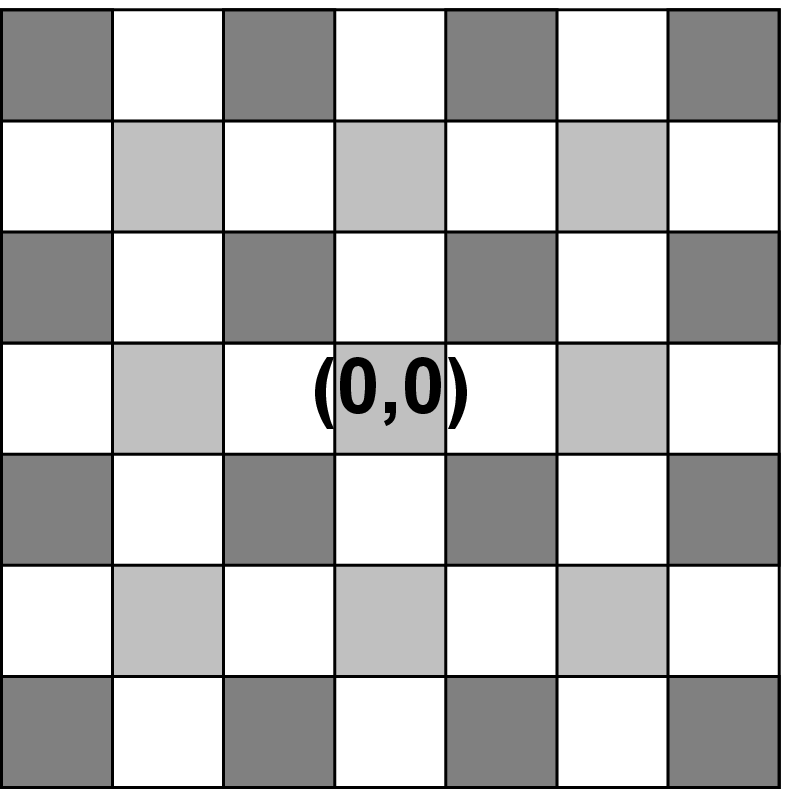} $\CA_1$}
\end{center}
\caption{Relation between ``coloring'' for $\CA_2$ and $\CA_1$} 
\label{fig:check2x}
\end{figure}
   
Relation between ``coloring'' for $\RCA_2$ ($\CA_2$) and $\RCA_1$ ($\CA_1$)
is shown on \Fig{check2x}. For $\RCA_2$ ($\CA_2$) coloring of cell $(i,j)$
used for illustration properties above is corresponding to value $i+j \mod 2$.
Next, all the $\CA_2$ board is mapped into sublattice of $\CA_1$
producing new type of coloring with {\em ``light''} and {\em ``dark''}
cells illustrated on \Fig{check2x}.

Due to such a map $\CA_2$ corresponds to sublattice $c'$ in $\CA_1$ with
coordinates $c_{i+j,i-j}$ \Eq{C2C1}. New indexes $(i+j,i-j)$
are both either odd or even. 

Let us use for $\RCA_1$ notation already introduced for $\RCA_2$ 
with $c_1$ and $c_2$ configurations corresponding to set of cells 
with nonzero first and second elements in the pair representing
a state of second-order CA. 

It was shown that for configurations derived from a single cell
with unit state such patterns have opposite color. For $\CA_2$
it corresponds to different diagonal sublattices and in each pattern
nonempty cells can not have adjoint sides, but may have common corners.
For $\CA_2$ with new scheme of coloring the corners of cells
are also separated.

\small

Let us first prove such expressions as \Eq{rec2R1} and \Eq{rec2R2}.
They already were derived above from \Eq{recR1} and \Eq{recR2}, 
but direct proof provided below illustrates some useful relations. 
The equation \Eq{rec2R2} may be derived from \Eq{TFT2n} and \Eq{FFT2n1}.
Let us recollect that for any polynomial $p(x,y)$ over GF(2) 
\begin{equation}
 p^2(x,y) = (\sum c_{i,j}x^i y^j)^2 = \sum (c_{i,j}x^i y^j)^2 =
 \sum (c_{i,j}x^{2i} y^{2j}) 
\label{p2xy} 
\end{equation}
and so for representations of two-states pattern via polynomials
used earlier the square corresponds to rescaling of the pattern 
$(i,j) \mapsto (2i,2j)$. The \Eq{TFT2n} corresponds to multiplication
of $T(x,y)$ on the rescaling pattern. For $\CA_1$ $T(x,y)$ 
is described by \Eq{TC1}. 

\begin{figure}[htb]
\begin{center}
\includegraphics[scale=0.4]{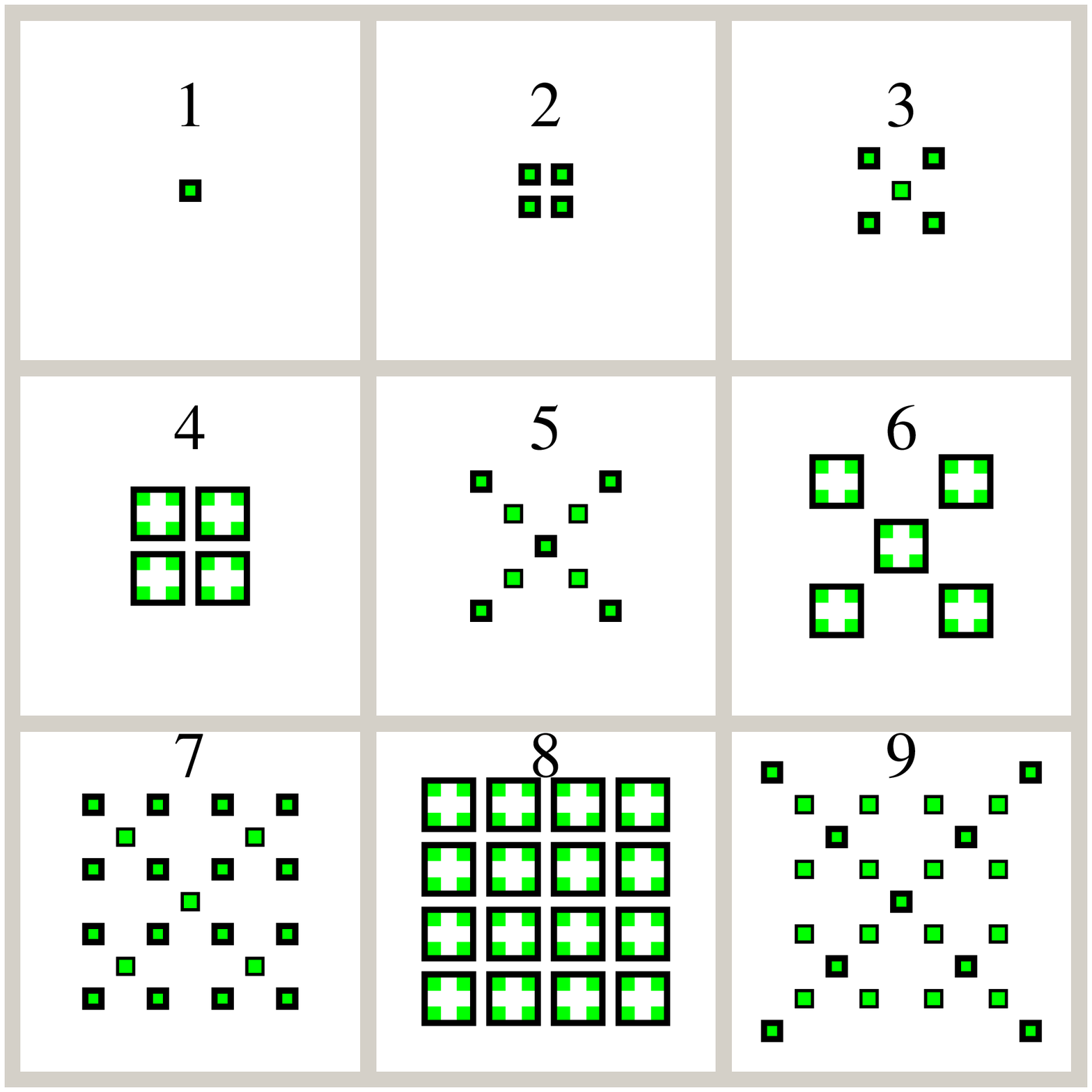}
\end{center}
\caption{Recursion \Eq{rec2R2} for $R_2(n)$ in $\CA_1$} 
\label{fig:compos2n}
\end{figure}

It was already shown, that for $\CA_1$ any cells with same value
are separated, so after the scaling distances between nonzero
cells are enough to put four new cells generated by $T(x,y)$
without overlap. 
\Fig{compos2n} illustrates that for 
$$n=1~{}_{\times 4} \to n=2~{}_{\times 4} \to n=4~{}_{\times 4} \to n=8,
\quad n=3~{}_{\times 4} \to n=6.$$ 
So, \Eq{TFT2n} proves first part of
\Eq{rec2R2}, $R_2(2n)=4R_2(n)$.

Next, due to \Eq{2nd} two polynomials $f_{n+1}$, $f_n$ in 
\Eq{FFT2n1} describes $(c_1, c_2)$ on a step $n$ and it was
already shown that the pattern are not intersecting for
chosen initial conditions. Square of the sum only rescales the
union without changing number of nonzero cells. 
\Fig{compos2n} illustrates that for 
$$n=1 \cup 2 \to n=3, ~ n=2 \cup 3 \to n=5, ~ n=3 \cup 4 \to n=7, 
~ n=4 \cup 5 \to n=9. $$ 
So,
\Eq{FFT2n1} proves second part of 
\Eq{rec2R2}, $R_2(2n+1)=R_2(n+1)+R_2(n)$. $\square$

\medskip

Recursive polynomial equation \Eq{f2kj} can be simply adopted for proof of
\Eq{recR2} for number of cells in $\RCA_1$ and $\RCA_2$ and it is enough
to demonstrate both \Eq{recR} and \Eq{recR1}.

Let us prove \Eq{recR2} for number of cells
with state 2 in $\RCA_1$ using \Eq{f2kj}.
The fact, that all cells with state 2 on each step $k$ 
are contained within a square region represented as direct
product of two {\em open} intervals $(-k,k) \times (-k,k)$
is also used and proved. 

For $k=0,1$ and initial configuration the \Eq{recR2} holds
and estimation for shape of square boundary is also true (for $k=0$
region is empty).

\begin{figure}[htb]
\begin{center}
\includegraphics[scale=0.4]{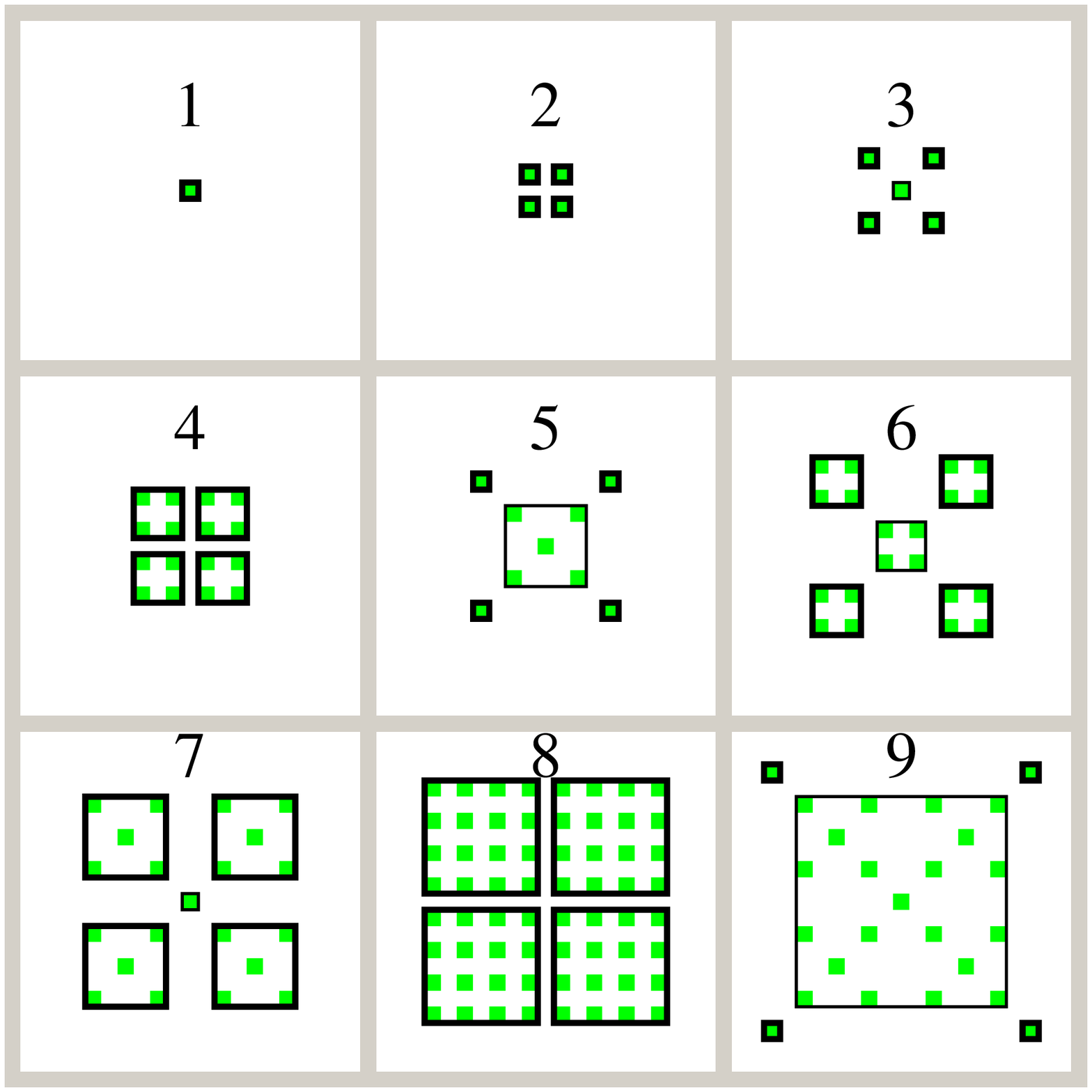}
\end{center}
\caption{Composition \Eq{recR2} for $R_2(n)$ in $\RCA_1$} 
\label{fig:compos2}
\end{figure}

Assume that equations hold for all patterns $j \le 2^n$ and consider
$j' = 2^n + j$. Due to \Eq{f2kj} and \Eq{C1expmov} the polynomial representation is
\begin{equation}
p_{2^n+j} =  
(x^{-2^n}y^{-2^n}\! +x^{2^n}y^{-2^n}\! + x^{-2^n}y^{2^n}\! +x^{2^n}y^{2^n})\,p_j
+ p_{2^n-j}.
\label{pC2nj}
\end{equation}
The multiplier before $p_j$ produces four copies moved in directions
$(-2^n,-2^n)$, $(-2^n,+2^n)$, $(+2^n,-2^n)$, $(+2^n,+2^n)$
and $p_{2^n-j}$ corresponds to pattern in the center, \Fig{compos2}. 
The five patterns are not overlapped: central one with 
$R_2(2^n-j)$ nonempty cells is contained
within $(-j'+j,j'-j) \times (-j'+j,j'-j)$ and other four others
with $R_2(j)$ nonempty cells are distributed within a ``four-fold'' 
disjointed region described by product
$$\{(-j',-j'+j)\cup(j'-j,j')\} \times \{(-j',-j'+j)\cup(j'-j,j')\}.$$
Total number of nonempty cells is $4 R_2(j) + R_2(2^n-j).$
So the equation for number of cells \Eq{recR2} holds for
$j' \le 2^{n+1}$. The union of the five regions
belongs to square $(-j',j') \times (-j',j')$.
$\square$

\Fig{compos2} illustrates relations
$$n=1~{}_{\times 4} \to n=2~{}_{\times 4} \to n=4~{}_{\times 4} \to n=8,$$ 
$$n=1~{}_{\times (4+1)} \to n=3, \quad n=2~{}_{\times (4+1)} \to n=6,$$
$$n=1_{\times 4}\cup 3 \to n=5, \quad n=3_{\times 4}\cup 1 \to n=7, 
\quad n=1_{\times 4}\cup 8 \to n=9. $$ 

\begin{figure}[htb]
\begin{center}
\includegraphics[scale=0.4]{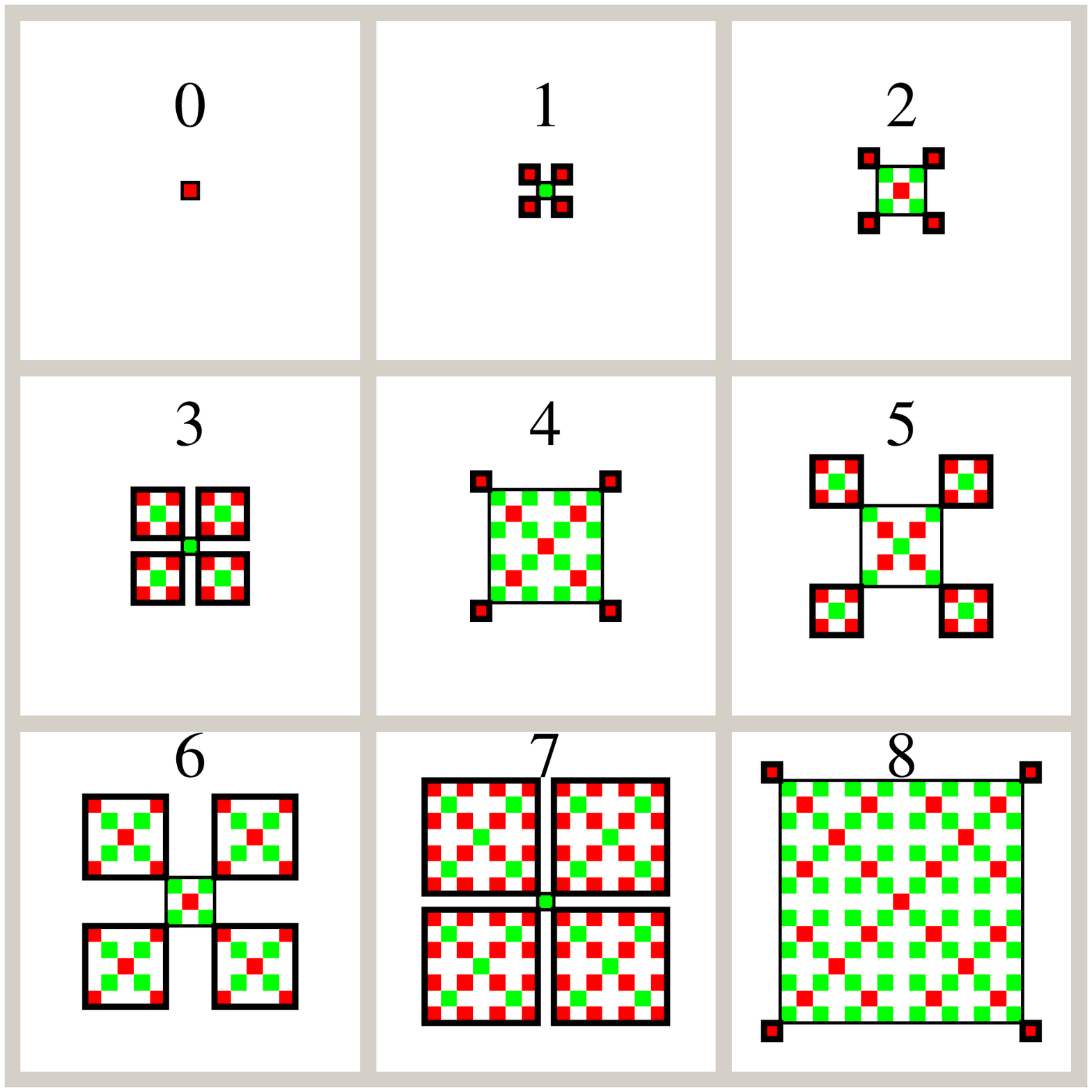}
\end{center}
\caption{Composition \Eq{recR} for $R(n)$ in $\RCA_1$} 
\label{fig:compos}
\end{figure}

The five patterns have a unit gap between them (\Fig{compos2}) and
only after consideration of all cells with nonzero values
corresponding to union of both ``checkerboard sublattices'' the
final patterns (\Fig{compos}) belong to square regions described
by product of {\em closed} intervals $[-k,k] \times [-k,k]$
and recursive equation \Eq{recR} corresponds to union
of five disjoint regions without gaps, \Fig{compos}.

\Fig{compos} illustrates relations
$$n=0_{\times 5} \to n=1,\quad n=0_{\times 4}\cup 1 \to n=2,\quad
n=0\cup 1_{\times 4}  \to n=3,$$
$$n=0_{\times 4}\cup 3 \to n=4,\quad 
n=1_{\times 4}\cup 2 \to n=5,\quad n=2_{\times 4}\cup 1 \to n=6, $$
$$n=0\cup 3_{\times 4} \to n=7,\quad n=0_{\times 4}\cup 7 \to n=8. $$

Let us check recursive equation for pair of polynomials \Eq{PFT}
representing all states of second-order CA and used for calculation of $R(n)$
\begin{eqnarray}
 P[C_{2^k+j}] &=& \bigl(f_{2^k+j+1}(T),f_{2^k+j}(T)\bigr) \notag \\
 &=& \bigl(T^{2^k}f_{j+1}(T)+f_{2^k-j-1}(T),T^{2^k}f_j(T)+f_{2^k-j}(T)\bigr) \notag \\
 &=& T^{2^k}\bigl(f_{j+1}(T),f_j(T)\bigr)+
     \bigl(f_{2^k-j-1}(T),f_{2^k-j}(T)\bigr) \notag \\
 &=& T^{2^k}  P[C_j] +  P[{\swap}C_{2^k-j-1}], 
\label{C2kj}    
\end{eqnarray}
where $\swap$ operation \Eq{swap} swaps values $1\leftrightarrow 2$. 

The \Eq{C2kj} illustrates dynamics of pattern growth, \Fig{compos}. 
Due to \Eq{invswap} application of transition rule $F$ to pattern
${\swap}C_i$ for any index $i > 0$ satisfies property
\begin{equation}
 F : \swap C_i \mapsto \swap C_{i-1},
\label{FXback} 
\end{equation}
so, application of $F$ to \Eq{C2kj} corresponds to increase of
four patterns $C_j$ and decrease of central region 
${\swap}C_{2^k-j-1}$ until $2^k-j-1 > 0$. 
For $j=2^k-1$ four outer configurations reach maximal size
and may not grow more, so on next step they are joined into
single central configurations $X C_{2^{k+1}-1}$ 
and four cells $C_0$ appear near corners as centers for future 
growth.

\begin{figure}[htb]
\begin{center}
\includegraphics[scale=0.4]{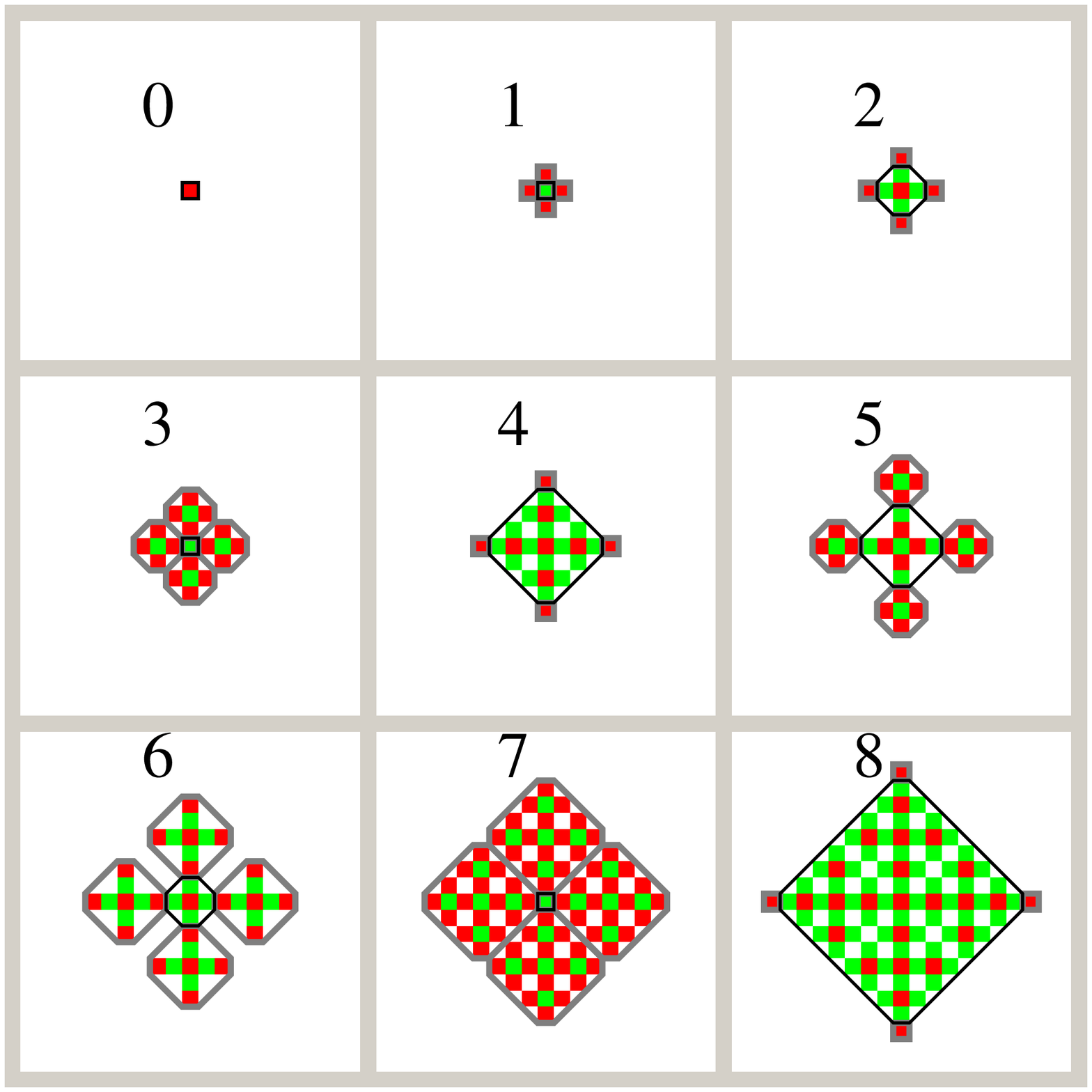}
\end{center}
\caption{Composition \Eq{recR} for $R(n)$ in $\RCA_2$} 
\label{fig:composR2}
\end{figure}

\medskip

Proofs of Eqs.~(\ref{recR}--\ref{rec2R2}) 
for $\RCA_2$ directly follow from consideration of $\RCA_1$, because
(similarly with relation between $\CA_1$ and $\CA_2$ discussed
earlier) $\RCA_2$ is equivalent with $\RCA_1$ acting on
{\em a diagonal sublattice}.
 
\begin{figure}[htb]
\begin{center}
\includegraphics[scale=0.4]{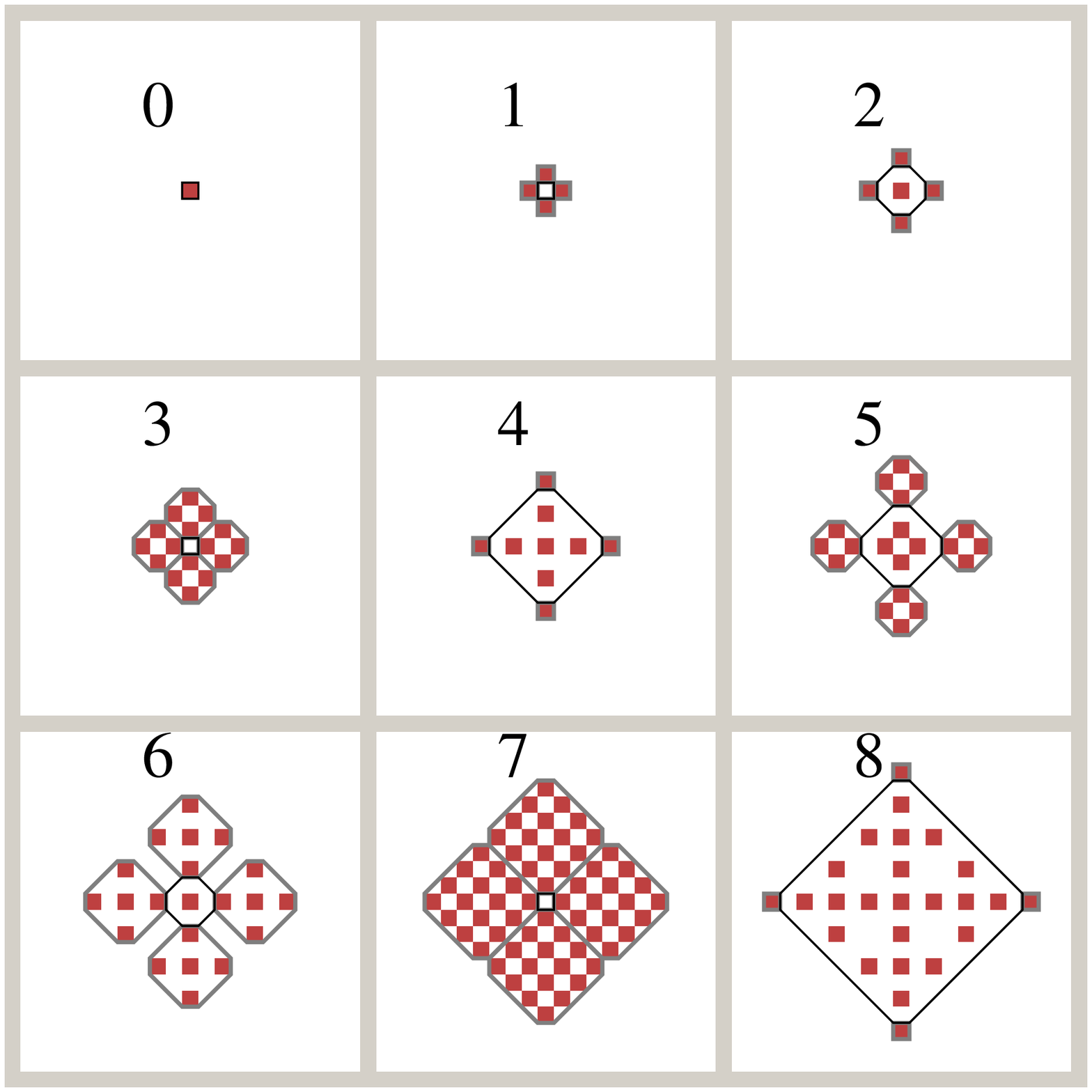}
\end{center}
\caption{$\RCA_2$, $\RCA_3$, $\RCA_3'$ --- cells with value 1} 
\label{fig:compos1R2}
\end{figure}

In such representation patterns
for $\RCA_2$ may look more closely packed \Fig{composR2},
but it does not change recursive equations due to above
mentioned equivalence. Let us now consider $\RCA_3$
and $\RCA_3'$.

Local rule for both $\RCA_2$ and $\RCA_3$ uses only four
closest cells with common sides in so-called {\em von Neumann
neighborhood}. Due to \Eq{2nd} it is enough to consider
actions of local rules for $\CA_2$ and $\CA_3$ on the
first element of pair to describe differences between rules. 
If the rules act in the same way for any configuration
under consideration, then actions of $\RCA_2$ and $\RCA_3$
for patterns derived from $\Delta_{0,0}$ are also the same.

Comparison of definition $\CA_2$ and $\CA_3$ shows
that local rules differ only for {\em three
nonempty cells} in von Neumann neighborhood.
On \Fig{compos1R2} for simplicity are shown only cells
with nonzero first components in the pair for configurations
used earlier, \Fig{composR2}.

All such pattern have 0,1,2,4 nonempty cells in von Neumann neighborhood
and so $\RCA_2$ and $\RCA_3$ act in the same way for such pattern.
Let us proof the property by induction. Any new configuration is 
composition of five previous patterns and it is enough to
consider new configurations near contiguities of they boundaries.

Due to consideration below for $n\neq 2^k-1$ there are four contacts
of central pattern with outer configurations. Four cells with two 
neighbors corresponds them. The cases $n = 2^k-1$ correspond
to contacts of four outer patterns and due to symmetry number
of neighbors there are always even. In fact, it may be simply shown
that all such configuration (of cells with state 1) are simple 
diamond-like checkerboard patterns with $2^k \times 2^k = 4^k$ 
cells, \Fig{compos1R2}.

Let us now consider $\RCA'_3$. The only difference between $\CA'_3$
and $\CA_3$ is additional requirement about cells with common 
corners. The limitation always holds due to ``coloring'' properties 
already discussed earlier on page \pageref{colprop}. 
Indeed, each new generation of cells with state 1 for $\RCA_2$ 
may appear only on checkerboard sublattice with opposite colors, 
{\em i.e.} all cells with common corner for an empty cell going
to be switched into the state 1 are empty.
So, evolution of $\RCA'_3$ starting with configuration $\Delta_{0,0}$
is also the same as for $\RCA_3$ and $\RCA_2$. $\square$

\end{document}